\documentstyle{amltd2004}
\begin{document}
\annalsline{155}{2002}
\received{August 8, 2000}
\startingpage{807}
\def\bye{\end{document}}
 \font\tenrm=cmr10

%--------------- Author macros ---------------
%for Bbb in amstex
\catcode`\@=11
\font\twelvemsb=msbm10 scaled 1100
\font\tenmsb=msbm10
%\font\ninemsb=msbm7 scaled 1100%msbm9
\font\ninemsb=msbm10 scaled 800
\newfam\msbfam
\textfont\msbfam=\twelvemsb  \scriptfont\msbfam=\ninemsb
  \scriptscriptfont\msbfam=\ninemsb
\def\msb@{\hexnumber@\msbfam}
\def\Bbb{\relax\ifmmode\let\next\Bbb@\else
 \def\next{\errmessage{Use \string\Bbb\space only in math
mode}}\fi\next}
\def\Bbb@#1{{\Bbb@@{#1}}}
\def\Bbb@@#1{\fam\msbfam#1}
\catcode`\@=12

 \catcode`\@=11
\font\twelveeuf=eufm10 scaled 1100
\font\teneuf=eufm10
\font\nineeuf=eufm7 scaled 1100%eufm9
\newfam\euffam
\textfont\euffam=\twelveeuf  \scriptfont\euffam=\teneuf
  \scriptscriptfont\euffam=\nineeuf
\def\euf@{\hexnumber@\euffam}
\def\frak{\relax\ifmmode\let\next\frak@\else
 \def\next{\errmessage{Use \string\frak\space only in math
mode}}\fi\next}
\def\frak@#1{{\frak@@{#1}}}
\def\frak@@#1{\fam\euffam#1}
\catcode`\@=12
%-------------- Author entries --------------------

%-------------- Article Text--------------------

%\intro %(Optional, Introduction)
\def\bt{\proclaim{Theorem}}
\def\et{\endproclaim}
\def\bl{\proclaim{Lemma}}
\def\el{\endproclaim}
\def\bp{\proclaim{Proposition}}
\def\ep{\endproclaim}

\newcommand{\A}{ {\cal A} }
\newcommand{\Ha}{ {\cal H} }
\newcommand{\Ls}{ {\cal L} }
\newcommand{\Es}{ {\cal E} }
\newcommand{\Fs}{ {\cal F} }
\newcommand{\Is}{ {\cal I} }
\newcommand{\Js}{ {\cal J} }
\newcommand{\Ks}{ {\cal K} }
\newcommand{\Cs}{ {\cal C} }
\newcommand{\Ms}{ {\cal M} }
\newcommand{\Pu}{ {\cal P} }
\newcommand{\Ts}{ {\cal T} }
\newcommand{\Se}{ {\cal S} }
\newcommand{\C}{{\Bbb C} }
\newcommand{\As}{{\Bbb A} }
\newcommand{\Hr}{{\Bbb H} }
\newcommand{\Ka}{{\Bbb K} }
\newcommand{\Oe}{{\Bbb O} }
\newcommand{\Ps}{{\Bbb P} }
\newcommand{\Q}{{\Bbb Q} }
\newcommand{\So}{{\Bbb S} }
\newcommand{\To}{{\Bbb T} }
\newcommand{\R}{{\Bbb R} }
\newcommand{\Xs}{{\Bbb X} }
\newcommand{\Z}{{\Bbb Z} }

\def\bx{{\bf x}}
\def\by{{\bf y}}
\def\bz{{\bf z}}
\def\bY{{\bf Y}}
\def\bG{{\bf G}}
\def\ba{{\bf a}}
\def\bb{{\bf b}}
\def\bh{{\bf h}}
\def\bce{{\bf c}}
\def\bte{{\bf t}}
\def\bu{{\bf u}}
\def\bv{{\bf v}}
\def\bw{{\bf w}}
\def\bal{\mathbold{\alpha}}
\def\bet{{\mathbold{\eta}}}
\def\bxi{{\mathbold{\xi}}}
\def\br{{\mathbold{\rho}}}
\def\bzeta{\mathbold{\zeta}}
\def\b0{{\bf 0}}

 \def\be{\begin{equation}}
\def\ee{\end{equation}}

\title{Linear equations in variables which lie\\ in a multiplicative group} %Article title
\shorttitle{Linear equations in a multiplicative group} % Shortened version for headline title

%% \acknowledgements{}
  \twoauthors{J.-H. Evertse, H. P. Schlickewei,}{W. M. Schmidt}

 \institutions{Universiteit te Leiden,  Leiden, The Netherlands\\
{\eightpoint {\it E-mail address\/}: evertse@math.leidenuniv.nl}\\
\vglue6pt
Universit\"at   Marburg, Marburg, Germany\\
{\eightpoint {\it E-mail address\/}: hps@mathematik.uni-marburg.de}\\
\vglue6pt
University of Colorado, Boulder, CO\\
{\eightpoint {\it E-mail address\/}:
 schmidt@euclid.colorado.edu}}

\centerline{\bf Abstract}
\vglue12pt

Let $K$ be a field of characteristic $0$ and let $n$ be a natural number. Let
$\Gamma$ be a subgroup of the multiplicative group $(K^*)^n$ of finite rank $r$.
Given $a_1,\ldots,a_n\in K^*$ write $A(a_1,\ldots,a_n,\Gamma)$ for the number of
solutions $\bx =\break (x_1,\ldots,x_n)\in \Gamma$ of the equation $a_1 x_1 +\cdots +
a_n x_n = 1$, such that no proper subsum of $a_1 x_1 +\cdots + a_n x_n$
vanishes. We derive an explicit upper bound for $A(a_1,\ldots,a_n,\Gamma)$ which
depends only on the dimension $n$ and on the rank $r$.

\section{Introduction} 

Let $K$ be an algebraically closed field of
characteristic $0$. Write $K^\ast$ for its multiplicative group of nonzero
elements, and let $(K^\ast)^n$ be the direct product consisting of $n$-tuples
$\bx=(x_1,\ldots,x_n)$ with $x_i\in K^\ast$ $\,(i=1,\ldots,n)$. So for
$\bx,\by\in (K^\ast)^n$ we write $\bx\ast\by=(x_1\,y_1,\ldots,x_n\,y_n)$.
Let $\Gamma$ be a subgroup of $(K^\ast)^n$ and suppose
$(a_1,\ldots, a_n)\in (K^\ast)^n$. We will be dealing with equations
\be
\label{1.1}
   a_1\, x_1 + \ldots + a_n\, x_n = 1
\ee
with $\bx\in \Gamma$.

A solution $\bx$ of (\ref{1.1}) is called {\it{nondegenerate}} if no subsum
of the left-hand side of (\ref{1.1}) vanishes, i.e., if $\sum\limits_{i\in I}
a_i\, x_i \neq 0$ for every nonempty subset $I$ of $\{1,\ldots, n\}$. Write
$A(a_1,\ldots, a_n; \Gamma)$ for the number of nondegenerate solutions
$\bx\in \Gamma$ of equation (\ref{1.1}).

Now suppose that $\Gamma$ has rank $r$. This means that there exists a finitely
generated subgroup $\Gamma_0$ of $\Gamma$, again of rank $r$, such that the
factor group $\Gamma/\Gamma_0$ is a torsion group. In other words, for any
$(x_1,\ldots,x_n)\in \Gamma$ there exists a natural number $k$ such that
$$
   \left(x_1^k,\ldots, x_n^k\right)\in\Gamma_0.
$$
We prove:

\bt
Suppose $\Gamma$ has finite rank $r${\rm .} Then the number $A(a_1,\ldots,\break a_n; 
\Gamma)$ of nondegenerate solutions $\bx\in \Gamma$ of equation
{\rm (\ref{1.1})}
satisfies the estimate
\be
\label{1.2}
   A(a_1,\ldots, a_n; \Gamma)\leq A(n,r) = \exp\left((6n)^{3n}(r+1)\right).
\ee
\et

The significant feature in our theorem is its uniformity. The bound (\ref{1.2})
depends only upon the dimension of the variety $V$ defined by equation
(\ref{1.1}) and upon the rank $r$ of the group $\Gamma$. We also remark that
once we have an estimate of the type
$$
   A(a_1,\ldots, a_n; \Gamma) \leq f(a_1,\ldots, a_n; n,r)
$$
with a function $f$ depending only on $a_1,\ldots,a_n$, $n$ and $r$, then we
get immediately
$$
   A(a_1,\ldots, a_n; \Gamma) \leq g(n, r)
$$
where $g$ is a function of $n$ and $r$ only. To see this, it suffices to
consider the equation
$$
   y_1+\cdots+y_n = 1
$$
and to ask for solutions $\by$ in the group generated by $(a_1,\ldots, a_n)$
and $\Gamma$ (which has rank $\leq r+1$).

It is conceivable that the function $A(n,r)$  we have given in (\ref{1.2}) is
far from best possible. In particular, no special care has been taken for the
numerical constants in (\ref{1.2}). However {\it{any function}}
$\widetilde A (n,r)$ which is suitable in (\ref{1.2}) indeed has to depend on
both $n$ and $r$.

As for the dependence on $n$ we give the following example. Pick elements 
$\alpha_1,\ldots,\alpha_n\in K^\ast$ with $\alpha_i\neq 1$ and $\alpha_i\neq 
\alpha_j$ $\,(1\leq i,\,j\leq n,\ i\neq j)$ and consider the equation
\be
\label{1.3}
   \left| 
   \begin{array}{lccclcr}
     1 & , & \ldots & , & 1 &,\ & 1\\
     \alpha_1 & , & \ldots & , & \alpha_n & ,\ & 1\\
     && \vdots &&&& \ \vdots\\
      \alpha_1^{n-1} & , & \ldots & , & \alpha_n^{n-1} & ,\ & 1\\
       x_1 & , & \ldots & , & x_n & ,\ & 1
     \end{array}\right| = 0.
\ee
This yields an equation
\be
\label{1.4}
   b_1\,x_1+\cdots+b_n\,x_n=1
\ee
with $b_i\in K^\ast$. But clearly (\ref{1.4}) has the $n$ solutions $\bx_i =
(x_{1i},\ldots,x_{ni})=(\alpha_1^i,\ldots,\alpha_n^i)$ $(i=0,\ldots,n-1)$.
Moreover, in the generic case, these will be nondegenerate solutions.
Therefore $\widetilde A(n,1)\geq n$. Bavencoffe and \pagebreak  B\'ezivin \cite{1} have
given a more sophisticated example which even shows that
$$
   \widetilde A(n,1)\geq c\,n^2
$$
where $c$ is an absolute constant.

On the other hand, suppose $n=p-1$ where $p$ is a prime. Let $\zeta$ be a
primitive $p$-th root of unity. Then
$$
   -\zeta-\zeta^2-\cdots-\zeta^{p-1}=1,
$$
and the same is true for any permutation of the roots on the left-hand side.
Therefore, for $n=p-1$ we have $\widetilde A(n,0)\geq n!$. We do not know
what should be in general the correct order of dependence on $n$ in
$\widetilde A(n,r)$.

As for the dependence on $r$, Erd\"{o}s, Stewart and Tijdeman \cite{7} have
constructed an example which shows that  
$$
   \widetilde A (2,r) \geq \exp \left( c\left( \frac{r}{\log\, 
   r}\right)^{\frac{1}{2}}\right)
$$
where $c$ is an absolute constant. This example may be extended to give
$$
  \widetilde A (n,r) \geq \exp \left( c (n)\left( \frac{r}{\log\, 
   r}\right)^{\frac{n-1}{n}}\right)
$$
where $c(n)$ depends only upon $n$. It has been conjectured that the correct
order of magnitude in $r$ should be of the shape
$$
   \exp \left( c (n)\left( \frac{r}{\log\, 
   r}\right)^{\frac{n}{n+1}}\right)
$$
or even
$$
   \exp \left( c (n)\ r^{\frac{n}{n+1}}\right).
$$

For $n = 2$, the assertion of Theorem 1.1 has been proved earlier.\break  Schlickewei
\cite{21} showed that $A(a_1,a_2,\Gamma)\leq c(r)$ and Beukers and Schlickewei
\cite{2} proved that we may take
\be
\label{1.5}
   c(r) = 2^{9(r+1)},
\ee
which clearly is much better than our bound $A(2,r)$.

For arbitrary $n$ and for $r=0$, i.e., when we are asking for solutions of
equation (\ref{1.1}) in roots of unity, Schlickewei \cite{20} proved that
we do not get more than $2^{4n!}$ nondegenerate solutions. This has been
considerably improved by Evertse \cite{11}. He obtained the bound
$$
   (n+1)^{3(n+1)^2},
$$
and this is much better than our bound $A(n,0)$ in (\ref{1.2}).\pagebreak

In all other cases, i.e., when $n\geq 3$ and $r\geq 1$, Theorem 1.1 is new. 
Previously, bounds involving only the dimension $n$ and the rank $r$ of the
group $\Gamma$ had been obtained only in the case when $\Gamma$ is the
$n$-fold product of
the group of $S$-units of a number field. We briefly review what was known in
the literature. Before we do so, let us remark that instead of the group
$\Gamma\subset 
(K^\ast)^n$ we could have considered a group $\Gamma'\subset K^\ast$ of finite
rank $r'$, say, and we could have asked for solutions of (\ref{1.1})
with $x_i\in \Gamma'$.
The difference is only minor, as the direct product $(\Gamma')^n$ then is a
subgroup of $(K^\ast)^n$ of rank $n\, r'$.

Writing $A'(a_1,\ldots, a_n; \Gamma')$ and $A'(n,r)$ for the quantities in
(\ref{1.2}) with respect to $\Gamma'\subset K^\ast$ we therefore see that
$$
   A'(n,r) \leq A(n, nr).
$$
The classical instances of equation (\ref{1.1}) are $S$-unit equations.
Let $F$ be a 
number field, let $S$ be a finite set of places of $F$ containing all the
archimedean ones and write $\Gamma(S)\subset F^\ast$ for the group of
$S$-units of $F$. For $n = 2$ 
and for $F = \Q$, Mahler \cite{14} has shown that
$$
  A'(a_1, a_2; \Gamma (S)) < \infty.
$$
Lang \cite{13} has extended Mahler's result to arbitrary number fields and
also to the case of arbitrary fields $K$ of characteristic $0$ and groups
$\Gamma\subset K^\ast$ of finite rank.

For general $n\geq 2$, Evertse \cite{9} and van der Poorten and Schlickewei
\cite{15} have shown that
$$
   A'(a_1,\ldots, a_n; \Gamma (S)) < \infty.
$$
The first quantitive result in our context is due to Evertse \cite{8}. He
proved for sets $S$ of cardinality $s$
\be
\label{1.6}
   A'(a_1, a_2;\Gamma (S)) \leq 3\cdot 7^{4s}.
\ee 
Notice that the group $\Gamma(S)$ is finitely generated and has rank $s-1$.
Therefore (\ref{1.6}) may be viewed as a special instance of a result of type
(\ref{1.2}) (cf. also (\ref{1.5})).

For arbitrary $n\geq 2$, Schlickewei \cite{17} proved that
\be
\label{1.7}
   A'(a_1,\ldots, a_n; \Gamma (S)) \leq c(n,s)
\ee
where $c$ is a function depending on $n$ and $s$ only. So (\ref{1.7}) again
is of the same type as (\ref{1.2}). The best explicit value for $c(n,s)$ is
due to Evertse \cite{10}. He proved
\be
\label{1.8}
   c(n,s) \leq 2^{35n^4 s}.
\ee
Now suppose $\Gamma$ is an arbitrary finitely generated subgroup of rank $r$ of
the multiplicative group $F^\ast$ of a number field $F$ of degree $d$. Taking
for $S$ the set of archimedean places of $F$ and those finite places whose associated
prime ideal divides some of the generators of $\Gamma$, we see that $\Gamma$
will be a subgroup of the group $\Gamma (S)$ of $S$-units. However the rank
$s-1$ of $\Gamma (S)$ may be much larger than the rank $r$ of the original
group $\Gamma$. So in general, even for groups $\Gamma\subset F^\ast$ the
bound $A'(n,r)$   we obtain with (\ref{1.2}) will be much better than the bound
of type (\ref{1.6}) or (\ref{1.8}) we get using the group $\Gamma (S)$.
Another disadvantage of $\Gamma(S)$ is the fact that $s=|S|\geq d/2$, $\,d$
being the degree of $F$. Therefore the device of estimating
$A'(a_1,\ldots,a_n;\Gamma)$ by $A'(a_1,\ldots,a_n;\Gamma(S))$   implicitly always
 introduces a dependence upon the degree of $F$ in the bound.

Schlickewei \cite{19} has estimated $A'(a_1,\ldots,a_n;\Gamma)$ in terms of
$n$, $r$ and $d$. And here Schlickewei and Schmidt \cite{24} have shown that
\be
\label{1.9}
   A'(a_n,\ldots, a_n; \Gamma) \leq (2d)^{41n^3r} r^{n^2r}.
\ee
The essential difference between (\ref{1.9}) and (\ref{1.2}) is the occurrence
of the degree $d$ in (\ref{1.9}). The problem in the current paper is to
estimate a quantity like the one on the left-hand side of (\ref{1.9})  avoiding
any dependence on $d$. We will come back to this at the end of this section.

It is well-known that results on equations (\ref{1.1}) are closely related to
results on multiplicities of linear recurrence sequences. A {\it{linear
recurrence sequence of order $n$}} is a sequence $\{u_m\}_{m\in\Z}$ of
elements in our field $K$ satisfying a relation
\be
\label{1.10}
   u_{m+n}=c_1\,u_{m+n-1}+\cdots+c_n\,u_m\quad (m\in\Z).
\ee
Here $c_1,\ldots,c_n$ are fixed elements from $K$. We assume that $n>0$ and
that relation (\ref{1.10}) is minimal, i.e., that $u_m$ does not satisfy a
relation of type (\ref{1.10}) for some $n'<n$. Then we have in particular
\be
\label{1.11}
            c_n\neq 0
\ee
(and $\{u_m\}$ is not the zero sequence). Define the
{\it{companion polynomial}} by
\be
\label{1.12}
    G(z)=z^n-c_1\,z^{n-1}-\cdots-c_n
        =\prod^r_{\rho=1}\,(z-\alpha_\rho)^{\sigma_\rho}
\ee
with distinct roots $\alpha_\rho$ of respective multiplicities $\sigma_\rho$
$\,(\rho=1,\ldots,r)$. By (\ref{1.11}), $\alpha_\rho\neq 0$ for
$\rho=1,\ldots,r$. Then we have a representation
\be
\label{1.13}
    u_m=\sum^r_{\rho=1} f_\rho(m)\,\alpha^m_\rho
\ee
where the $f_\rho$ are polynomials. It follows from the minimality of relation 
(\ref{1.10}) that $f_\rho(x)$ has degree $\sigma_\rho-1$ $\,(\rho=1,\ldots,r)$. 
The sequence $\{u_m\}$ is called {\it{nondegenerate}} if no quotient
$\alpha_i/\alpha_j$ $\,(1\leq i<j\leq r)$ is a root of unity.

We say that the sequence $\{u_m\}$ is {\it{simple}} if the companion polynomial
$G(z)$ has only simple zeros. In that case the quantities $\sigma_\rho$ in
(\ref{1.12}) are \pagebreak all equal to $1$, so the polynomials $f_\rho$ in (\ref{1.13})
are constants and we have
\be
\label{1.14}
    u_m = a_1\,\alpha_1^m+\cdots+a_n\,\alpha_n^m\quad(m\in\Z)
\ee
with nonzero coefficients $a_i\in K$ and with distinct elements
$\alpha_i\in K^\ast$.

Write ${\cal S}(u_m)$ for the set of zeros of $\{u_m\}$, i.e., for the set of
solutions $k\in\Z$ of the equation
\be
\label{1.15}
    u_k = 0.
\ee
When $\{u_m\}$ has order $1$, then trivially ${\cal S}(u_m) = \emptyset$.
Therefore from now on we will only consider sequences $\{u_m\}$ of order
$n\geq 2$.

The classical theorem of Skolem-Mahler-Lech says that for arbitrary linear
recurrence sequences $\{u_m\}$ of order $\geq 2$, $\,{\cal S}(u_m)$ is the
union of a finite set of integers and a finite number of arithmetic
progressions. This implies in particular that for nondegenerate sequences
$\{u_m\}$ the set ${\cal S}(u_m)$ is finite.

An old conjecture says that for nondegenerate sequences $\{u_m\}$ of order
$n\geq 2$ the cardinality of ${\cal S}(u_m)$ is bounded in terms of $n$ only.
For $n=2$, by nondegeneracy it is obvious that $|{\cal S}(u_m)|\leq 1$.
Schlickewei \cite{22} proved the conjecture for $n=3$. Beukers and Schlickewei
\cite{2} derived for nondegenerate sequences $\{u_m\}$ of order $3$ the bound
$$
    |{\cal S}(u_m)|\leq 61.
$$
For nondegenerate sequences $\{u_m\}$ of {\it{rational numbers}} and of
arbitrary order~$n$, the conjecture was proved by Schlickewei \cite{18}.

We now study {\it{simple}} recurrence sequences $\{u_m\}$ (never mind whether
degenerate or not). For such sequences, in view of (\ref{1.14}), equation
(\ref{1.15}) becomes
\be
\label{1.16}
    a_1\,\alpha_1^k+\cdots+a_n\,\alpha_n^k = 0\quad (k\in\Z).
\ee
Applying Theorem 1.1 to groups $\Gamma$ of rank $\leq 1$ we deduce:
\bt \hskip-8pt
Let $K$ be an algebraically closed field of characteristic~$0${\rm .} Suppose
$n\geq 3$ and let $\{u_m\}_{m\in\Z}$ be a simple linear recurrence sequence in
$K$ of order $n${\rm .} Then there are integers $k_1,\ldots,k_{q_1}$ and arithmetic
progressions $T_1,\ldots,T_{q_2}$ of the shape
$$
   T_i = \{a_i+t\,v_i\ |\ t\in\Z\},\quad a_i,v_i\in\Z,\quad v_i\neq 0\quad
   (i=1,\ldots,q_2),
$$
where
\be
\label{1.17}
   q_1+q_2\leq\exp\left((6n)^{3n}\right),
\ee
such that
$$
   {\cal S}(u_m) = \{k\in\Z\ \,|\ u_k=0\} = \{k_1,\ldots,k_{q_1}\}\cup T_1
\cup\ldots\cup T_{q_2}.
$$
In particular{\rm ,} if $\{u_m\}$ is nondegenerate{\rm ,} then ${\cal S}(u_m)$ has
cardinality
\be
\label{1.18}
    |{\cal S}(u_m)|\leq\exp\left((6n)^{3n}\right).
\ee
\et
\vglue-6pt
Theorem 1.2 is a uniform quantitative version of the Skolem-Mahler-Lech
theorem. In the meantime, W. M. Schmidt \cite{27} has proved that {\it{for any
nondegenerate sequence}} $\{u_m\}$ (\/{\it{even if not simple\/{\rm )} the set}} 
${\cal S}(u_m)$ {\it{has cardinality bounded in terms of the order}} $n$ 
{\it{only}}.

The bound obtained by Schmidt in this more general setting is triply
exponential in $n$. Moreover, in his recent paper \cite{28}, Schmidt has
also proved that Theorem 1.2 is true in general and not only for simple
sequences. However, again instead of (\ref{1.17}) he gets a bound which is
triply exponential in $n$.

The new ingredients in our proof are as follows. On the one hand we apply the
absolute version of the Subspace Theorem due to Evertse and Schlickewei \cite{12}. 
On the other hand we use a result of Schmidt \cite{26} on lower bounds for heights of
points on varieties. 

In proving our theorems, by a specialization argument we may restrict
ourselves to the situation when in (\ref{1.1}) (or in (\ref{1.16})
respectively) all quantities involved are algebraic. Indeed it suffices to
prove a result of type (\ref{1.9}), but without dependence upon the degree
$d$ of the number field.

An application of the Subspace Theorem then gives an assertion on the ``large'' 
solutions of equation (\ref{1.1}). In fact the bound it gives for the number of 
``large'' solutions involves only   the ``good'' parameters $n$ and $r$. So 
for the quantitative result all depends upon the parameters showing up in the 
definition of ``small''. Usually in this definition the parameters $n$ and $d$ 
showed up. Thus in estimating the number of ``small'' solutions the parameter 
$d$ could not be avoided. In a recent paper \cite{12}, Evertse and Schlickewei 
have proved a new absolute quantitative version of the Subspace Theorem which 
in turn makes use of the absolute Minkowski Theorem established by Roy and 
Thunder \cite{16}. The definition of ``small'' in the absolute Subspace Theorem 
does not depend upon the degree $d$ at all.

Unfortunately this does not suffice yet. To handle the ``small'' solutions
usually one applies a gap principle. For this purpose one needs a {\it{lower}}
bound for the height of a small solution. Traditionally, this was achieved via
Dobrowolski's theorem \cite{6}. But here again the degree $d$ comes in. To
overcome this difficulty we apply lower bounds for heights of points on
varieties as given in recent work of Zhang \cite{29}, Bombieri and Zannier
\cite{3}, and in explicit form for the first time by Schmidt \cite{26}.

\vglue4pt {\it {R}emark}.  In recent work \cite{4}, \cite{5}, David and Philippon
  proved a slight sharpening of Schmidt's results \cite{26}. It is easily
seen that with this sharpening the bound for $A(n,r)$ given in (\ref{1.2})
can be improved to
$$
     A(n,r)\leq\exp\left((r+1)\exp(c_1n)\right).
$$
Similarly, the bound (\ref{1.17}) can be improved to
$$
    q_1+q_2\leq\exp\exp(c_1n).
$$
Here $c_1$ is an absolute constant.
 
\section{Algebraic points}

In the case when in (\ref{1.1}) all quantities involved are algebraic we can
prove a slightly more general result.

Let $F$ be a number field. Write $M(F)$ for the set of its places. For each
$v\in M(F)$ we let $|\ \,|_v$ be the associated absolute value such that for
$x\in\Q$ we have
\be
\label{2.1}
  |x|_v = \left\{
             \begin{array}{ll}
                |x| & \mbox{if}\ \ v|\infty\\
                |x|_p & \mbox{if}\ \ v|p\,,
             \end{array}
           \right.
\ee
where $p$ is a prime number and where $|p|_p = p^{-1}$. We denote the
completion of $F$ at the place $v$ by $F_v$; similarly for $p\in M(\Q)$,
$\,\Q_p$ denotes the completion of $\Q$ at $p$ (so that $\Q_{\infty}=\R$, the
field of real numbers). The {\it{normalized}} absolute value $\|\ \,\|_v$ on
$F$ then is defined by
\be
\label{2.2}
    \|x\|_v={|x|_v}^{[F_v:\Q_p]/[F:\Q]}\quad\mbox{if}\ \ v|p.
\ee
We write $\overline{\Q}$ for the algebraic closure of $\Q$. Given $\bx=
(x_1,\ldots,x_n)\in\overline{\Q}^{\,n}$, we define the absolute multiplicative
height $H(\bx)$ as follows: we choose a number field $F$ such that $\bx\in
F^n$ and we put
\be
\label{2.3}
    H(\bx)=\prod_{v\in M(F)}\max\{1,\|x_1\|_v,\ldots,\|x_n\|_v\}.
\ee
Notice that (\ref{2.3}) does not depend on the choice of $F$. We define the
absolute logarithmic height $h(\bx)$ by
\be
\label{2.4}
    h(\bx)=\log H(\bx).
\ee

In \cite{23}, Schlickewei and Schmidt proved the following result.

 {\it{Let
$F$ be a number field of degree $d$. Let}}
\be
\label{2.5}
    \Gamma\subset(F^*)^n\ \mbox{\it{be a finitely generated subgroup with}}
        \ \,\mbox{rank}\,\Gamma=r. \hskip.4in
\ee
{\it{Consider the equation}}
\be
\label{2.6}
    y_1+\cdots+y_n=1,
\ee
{\it{to be solved in vectors $\by=(y_1,\ldots,y_n)\in F^n$ of the shape}}
\be
\label{2.7}
 \by=\bx*\bz\ \ \mbox{\it{with}}\ \ \bx\in\Gamma,\ \,\bz\in(\Q^*)^n,\ \,h(\bz)
     \leq\frac1{4n^2}\,h(\bx).
\ee
{\it{Then the set of solutions $\by$ of}} (\ref{2.6}), (\ref{2.7}) {\it{is
contained in the union of not more than}}
\be
\label{2.8}
    2^{30n^2}\left(32n^2\right)^rd^{3r+2n}
\ee
{\it{proper linear subspaces of}} $F^n$.

Instead of (\ref{2.5}), we now suppose
\be
\label{2.9}
    \Gamma\ \mbox{is a subgroup of}\ {\left(\,\overline{\Q}^{\,*}\right)}^n
    \ \mbox{of rank}\ r.
\ee
So now $\Gamma$ is not necessarily finitely generated. (On the other hand, we
notice that (\ref{2.9}) is more special than the setting studied in Section 1,
where we assumed $\Gamma\subset(K^*)^n$ for some algebraically closed field
$K$ of characteristic $0$, so that in fact implicitly we assumed that
$\overline{\Q}\subset K$.)

Again we consider equation (\ref{2.6}). However, instead of (\ref{2.7}) we now
ask for solutions $\by\in{(\,\overline{\Q}^{\,*})}^n$ of the shape
\be
\label{2.10}
    \by=\bx*\bz\ \ \mbox{with}\ \ \bx\in\Gamma,
    \ \,\bz\in{
\left(\,\overline{\Q}^{\,*}\right)}^n,
    \ \,h(\bz)\leq n^{-1}\exp\left(-(4n)^{3n}\right)(1+h(\bx) ).
\ee
We prove:

\bt
Let $n\geq 2${\rm .} Suppose that $\Gamma$ is a subgroup of
${(\,\overline{\Q}^{\,*})}^n$ of finite rank $r${\rm .} Then the set of points
$\by\in\overline{\Q}^{\,n}$ satisfying {\rm (\ref{2.6})} and {\rm (\ref{2.10})}
is contained in the union of not more than
\be
\label{2.11}
    B(n,r)=\exp\left((5n)^{3n}(r+1)\right)
\ee
proper linear subspaces of $\overline{\Q}^{\,n}${\rm .}
\et

It turns out that Theorem 1.1 as well as Theorem 1.2 follow from Theorem~2.1.
Indeed in Section 3 we give a specialization argument which reduces the
situation we encounter in Section 1 to a setting where all quantities are
algebraic. In Section 4 we then prove Theorem 1.1 by means of induction using
Theorem 2.1. In Section 5, Theorem 1.2 will be deduced from Theorem 1.1. The
remainder of the paper, starting with Section 6, then is devoted to the proof
of Theorem 2.1.
\vglue-8pt
\section{Specialization}

\vglue-4pt
Let $K$ be the field from Section 1. Since $K$ is algebraically closed and has
characteristic equal to zero, we may suppose that $\overline{\Q}\subset K$.

\bl
Let $U=\{u_1,\ldots,u_k\}$ be a finite subset of $K${\rm .} Then there exists a
ring homomorphism
\be
\label{3.1}
    \varphi\,:\,\overline{\Q}\,[U]\,\longrightarrow\,\overline{\Q}
\ee
whose restriction to $\overline{\Q}$ is the identity{\rm .}
\el

{\it Proof}.  We recall the proof of this well-known fact. Let
$\Js$ be the ideal of polynomials $f\in\overline{\Q}\,[X_1,\ldots,X_k]$
with $f(u_1,\ldots,u_k)=0$. Clearly $1\notin\Js$ and therefore $\Js\neq
\overline{\Q}\,[X_1,\ldots,X_k]$. Thus by Hilbert's Nullstellensatz there
exists a point $\bce=(c_1,\ldots,c_k)\in\overline{\Q}^{\,k}$ with $f(\bce)=0$
for each $f\in\Js$. The ring $\overline{\Q}\,[U]=\overline{\Q}
\,[u_1,\ldots,u_k]$ consists of all expressions $g(u_1,\ldots,u_k)$ with
$g\in\overline{\Q}\,[X_1,\ldots,X_k]$. We consider the diagram
$$
\overline{\Q}\,[u_1,\ldots,u_k]\,\longrightarrow\,\overline{\Q}\,[X_1,\ldots,
X_k]/\Js\,\longrightarrow\,\overline{\Q}
$$
where the mappings are given by
$$
g(u_1,\ldots,u_k)\,\longmapsto\,g\ \mbox{mod}\ \Js\,\longmapsto
\,g(c_1,\ldots,c_k).
$$
These mappings are well-defined ring homomorphisms leaving $\overline{\Q}$
invariant. Their composition yields the desired homomorphism $\varphi$ in
(\ref{3.1}).

In order to prove Theorem 1.1, it will suffice to show that any {\it{finite}}
subset $M$ of the set of nondegenerate solutions of equation (\ref{1.1}) has
cardinality
\be
\label{3.2}
    \leq A(n,r).
\ee
Write $M=\{\bx_1,\ldots,\bx_m\}$ with $\bx_i=(x_{i1},\ldots,x_{in})$ $\,(i=
1,\ldots,m)$. We want to map $M$ injectively to a set of nondegenerate
solutions of an equation of type (\ref{1.1}) where, however, $\ba=(a_1,\ldots,
a_n)\in\overline{\Q}^{\,n}$ and $\Gamma\subset{(\,\overline{\Q}^{\,*})}^n$.
We will then be in a position to apply Theorem~2.1.

Let $U=\{u_1,\ldots,u_k\}\subset K$ be the set consisting of the following
elements:
\be
\label{3.3}
    a_1,\ldots,a_n;
\ee
\be
\label{3.4}
    x_{ij}\quad (i=1,\ldots,m;\,j=1,\ldots,n);
\ee
\be
\label{3.5}
    \sum_{j\in I}a_j\,x_{ij}\quad (i=1,\ldots,m;\,I\subset\{1,\ldots,n\},\,
        I\neq\emptyset);
\ee
\be
\label{3.6}
    x_{i_1,j}-x_{i_2,j}\quad (1\leq i_1<i_2\leq m;\,j=1,\ldots,n);
\ee
\be
\label{3.7}
    \mbox{the multiplicative inverses of all nonzero numbers in (\ref{3.3})--(\ref{3.6}).}\hskip.25in
\ee
Let $\varphi$ be a ring homomorphism from $\overline{\Q}\,[U]$ into $\overline
{\Q}$ as in Lemma 3.1. By (\ref{3.7}), the nonzero elements in $U$ are units
in the ring $\overline{\Q}\,[U]$. Therefore they are mapped by $\varphi$ to
nonzero elements of $\overline{\Q}$.

Write $a'_j=\varphi(a_j)$, $x'_{ij}=\varphi(x_{ij})$, $\bx'_i=(x'_{i1},\ldots,
x'_{in})$ $\,(i=1,\ldots,m;\,j=1,\ldots,n)$. Then by (\ref{1.1}) we get
\be
\label{3.8}
  a'_1\,x'_{i1}+\cdots+a'_n\,x'_{in}=\varphi\left(\sum_{j=1}^na_j\,x_j\right)
         =1\quad (i=1,\ldots,m).
\ee
The numbers in (\ref{3.5}), by nondegeneracy, are nonzero. Therefore their
images under $\varphi$ are nonzero as well. We may conclude that
\be
\label{3.9}
 \sum_{j\in I}a'_j\,x'_{ij}\neq 0\quad (i=1,\ldots,m;\,I\subset\{1,\ldots,n\},
          \,I\neq\emptyset).
\ee
Moreover, the nonzero numbers in (\ref{3.6}) have nonzero images. This
implies that $\bx'_1,\ldots,\bx'_m$ are distinct.

Let $\Gamma_1$ be the subgroup of $\Gamma$ generated by $\bx_1,\ldots,\bx_m$.
Then $\Gamma_1$ has rank $\leq r$. We infer from (\ref{3.7}) that $\Gamma_1
\subset{(\,\overline{\Q}\,[U])}^n$. Let $\Gamma'_1$ be the multiplicative
subgroup of ${(\,\overline{\Q}^{\,*})}^n$ generated by $\bx'_1,\ldots,\bx'_m$.
Then $\Gamma'_1$ is the image of $\Gamma_1$ under the group homomorphism
$$
(x_1,\ldots,x_n)\,\longmapsto\,(\varphi(x_1),\ldots,\varphi(x_n)).
$$
We may conclude that $\Gamma'_1$ has rank $\leq r$.

Altogether we see that $\bx'_1,\ldots,\bx'_m$ are distinct, nondegenerate
solutions of the equation
\be
\label{3.10}
   a'_1\,x'_1+\cdots+a'_n\,x'_n=1
\ee
to be solved in vectors
\be
\label{3.11}
   \bx'=(x'_1,\ldots,x'_n)\in\Gamma'_1.
\ee
Here $a'_1,\ldots,a'_n\in\overline{\Q}^{\,*}$ and $\Gamma'_1$ is a subgroup of
${(\,\overline{\Q}^{\,*})}^n$ of rank $\leq r$.

Notice that $A(n,r)$ in (\ref{1.2}) satisfies
$$
   A(n,r_1)<A(n,r_2)\ \ \mbox{for}\ \ r_1<r_2.
$$
Therefore we have shown:
\bl
In order to prove Theorem {\rm  1.1,} we may   suppose without loss of
generality  that
$K=\overline{\Q}${\rm .}
\el

\section{Deduction of Theorem 1.1 from Theorem 2.1}

In view of Lemma 3.2, we may suppose that $K=\overline{\Q}$. Under this
hypothesis we show that equation (\ref{1.1}) does not have more than
$$
    A(n,r)
$$
nondegenerate solutions $\bx\in\Gamma$, where
\be
\label{4.1}
    A(n,r)=\exp\left((6n)^{3n}(r+1)\right)
\ee
as in (\ref{1.2}).

The case $n=1$ is obvious. Now suppose $n>1$ and our claim  to be shown for
$n'<n$.

Let $B(n,r)$ be the quantity from (\ref{2.11}) in Theorem 2.1. Write $\Gamma'$
for the group generated by $\ba=(a_1,\ldots,a_n)$ and $\Gamma$. So if $\bx$
runs through $\Gamma$, the point $\by=\ba*\bx$ runs through $\Gamma'$. Clearly
$\Gamma'$ has rank $\leq r+1$. Thus the solutions $\bx\in\Gamma$ of (\ref{1.1})
give rise to solutions $\by\in\Gamma'$ of the equation
\be
\label{4.2}
    y_1+\cdots+y_n=1.
\ee
Applying Theorem 2.1 with $\bz = (1,\ldots,1)$ to equation (\ref{4.2}) and the 
group $\Gamma'$, we may
infer that the set of solutions $\by\in\Gamma'$ of (\ref{4.2}) (never mind
whether degenerate or not) is contained in the union of $B(n,r+1)$ proper
linear subspaces of $\overline{\Q}^{\,n}$. Consequently, also the set of
solutions $\bx\in\Gamma$ of equation (\ref{1.1}) is contained in the union of not more than
\be
\label{4.3}
    B(n,r+1)
\ee
proper linear subspaces of $\overline{\Q}^{\,n}$.

Let $V$ be one of these subspaces, defined by an equation
\be
\label{4.4}
    \sum_{i\in I}b_i\,x_i=0
\ee
where $I$ is a subset of $\{1,\ldots, n\}$ of cardinality $|I|\geq 2$, and
where $b_i\neq 0$ for $i\in I$. Let $J$ be a nonempty subset of $I$ and
consider those $\bx\in\Gamma\cap V$ for which
\be
\label{4.5}
   \sum_{i\in J}b_i\,x_i=0,
\ee
but no proper nonempty subsum of (\ref{4.5}) vanishes. Thus $2\leq|J|\leq n$.

Let us suppose for the moment that $J=\{1,\ldots,\ell\}$. Writing
$c_i=-b_i/b_1$ we get, with $w_i=x_i/x_1$ $\,(i=2,\ldots,\ell)$,
\be
\label{4.6}
   \sum^\ell_{i=2}c_i\,w_i=1.
\ee
Now $(x_1,\ldots,x_\ell,x_{\ell+1},\ldots,x_n)\in\Gamma$; therefore
$(w_2,\ldots,w_\ell)$ lies in the group $\Gamma_1$ consisting of
$(\ell-1)$-tuples such that
\be
\label{4.7}
   (u,u\,w_2,\ldots,u\,w_\ell,u_{\ell+1},\ldots,u_n)\in\Gamma
\ee
for some $u,u_{\ell+1},\ldots,u_n$. Let $\Gamma_2$ be the group of elements
\be
\label{4.8}
   (x,\ldots,x,x_{\ell+1},\ldots,x_n)\in\Gamma.
\ee
The map
$$
   (x_1,x_2,\ldots,x_\ell,x_{\ell+1},\ldots,x_n)\,\longmapsto\,\left(
          \frac{x_2}{x_1},\ldots,\frac{x_\ell}{x_1}\right)
$$
is a surjective homomorphism $\Gamma\to\Gamma_1$ with kernel $\Gamma_2$.
Therefore, when $\mbox{rank}\,\Gamma_i=r_i$ we have $r_1+r_2=r$.

By induction the equation (\ref{4.6}) has at most $A(\ell-1,r_1)$
nondegenerate solutions. When $(w_2,\ldots,w_\ell)$ is such a solution, fix
$u,u_{\ell+1},\ldots,u_n$ with (\ref{4.7}). The original solution $\bx$ of
(\ref{1.1}) is of the form
\be
\label{4.9}
    (x,x\,w_2,\ldots,x\,w_\ell,x_{\ell+1},\ldots,x_n),
\ee
so that
\be
\label{4.10}
   b\,x+\sum^n_{i=\ell+1}a_i\,x_i=1
\ee
with $b=a_1+\sum\limits^\ell_{i=2}a_i\,w_i$. If the solution $\bx$ of
(\ref{1.1}) is nondegenerate, then so is the solution $(x,x_{\ell+1},\ldots,
x_n)$ of (\ref{4.10}).

Taking the quotient of (\ref{4.7}), (\ref{4.9}) we see that $$(x/u,\ldots,x/u,
x_{\ell+1}/u_{\ell+1},\ldots,x_n/u_n)\in\Gamma_2.$$ With the notation $x'=x/u$,
$\,x'_i=x_i/u_i$ $\,(i=\ell+1,\ldots,n)$, (\ref{4.10}) becomes
\be
\label{4.11}
    b'x'+\sum^n_{i=\ell+1}a'_i\,x'_i=1
\ee
where $b'=b\,u$, $\,a'_i=a_i\,u_i$ $\,(i=\ell+1,\ldots,n)$. By induction,
and since $n-\ell+1\break <n$, (\ref{4.11}) has not more than $A(n-\ell+1,\,r_2)$
nondegenerate solutions. Combining this with the bound $A(\ell-1,\,r_1)$ for
the number of solutions of (\ref{4.6}), we see that (\ref{4.5}) gives rise to
not more than
\be
\label{4.12}
   A(\ell-1,\,r_1)\,A(n-\ell+1,\,r_2)\leq A(n-1,\,r)
\ee
solutions of (\ref{1.1}); the last inequality is a consequence of 
$$
   A(a,r_1)\,A(b,r_2)\leq A(a+b-1,\,r_1+r_2),
$$
which follows from the definition (\ref{4.1}) of $A(n,r)$. Taking account of
the possible subsets $J$ of $I$, we see that each subspace $V$ contains at
most $2^nA(n-1,\,r)$ solutions. We still have to multiply this by the number
$B(n,\,r+1)$ of subspaces. In this way we obtain a bound
$$
   2^nA(n-1,\,r)B(n,\,r+1).
$$
This is
$$
   2^n\exp\left((6(n-1))^{3(n-1)}(r+1)\right)\exp\left((5n)^{3n}(r+2)\right)<\exp\left((6n)^{3n}(r+1)\right),
$$
and Theorem 1.1 follows.\pagebreak

\section{Proof of Theorem 1.2}

Let $\{u_m\}$ be a simple linear recurrence sequence of order $n\geq 2$
contained in an algebraically closed field $K$ of characteristic $0$. To
simplify our exposition, single elements of $\Z$ will also be called
arithmetic progressions (indeed they may be viewed as arithmetic progressions
with difference $0$). Thus we have to show that the set
$$
    {\cal S}(u_m)=\{k\in\Z\ \,|\ u_k=0\}
$$
is the union of at most
\be
\label{5.1}
    W(n)=\exp\left((6n)^{3n}\right)
\ee
arithmetic progressions.

We proceed by induction on $n$. For $n=2$, our assertion is obvious. Assume
$n\geq 3$. Recall that
$$
   u_m=a_1\,\alpha_1^m+\cdots+a_n\,\alpha_n^m
$$
for certain nonzero elements $a_1,\ldots,a_n,\alpha_1,\ldots,\alpha_n\in K$.
Hence ${\cal S}(u_m)$ is the set of solutions $k\in\Z$ of
\be
\label{5.2}
    a_1\,\alpha_1^k+\cdots+a_n\,\alpha_n^k=0.
\ee
First consider those $k\in\Z$ for which no proper subsum of the left-hand side
of (\ref{5.2}) vanishes. For each such $k$, the vector
$$
   \left((\alpha_1/\alpha_n)^k,\ldots,(\alpha_{n-1}/\alpha_n)^k\right)
$$
is a nondegenerate solution of
\be
\label{5.3}
    \left(-\frac{a_1}{a_n}\right)x_1+\cdots+\left(-\frac{a_{n-1}}{a_n}\right)
        x_{n-1}=1\ \ \mbox{in}\ \ \bx=(x_1,\ldots,x_{n-1})\in\Gamma, \qquad
\ee
where $\Gamma$ is the group generated by $(\alpha_1/\alpha_n,\ldots,
\alpha_{n-1}/\alpha_n)$. Clearly $\Gamma$ has rank $\leq 1$. So by Theorem 1.1,
equation (\ref{5.3}) has at most
\be
\label{5.4}
    A(n-1,1)=\exp\left((6(n-1))^{3(n-1)}2\right)
\ee
nondegenerate solutions. As can be easily verified, for each solution\break
$(x_1,\ldots,x_{n-1})$ of (\ref{5.3}) the set of $k\in\Z$ with $\left(
(\alpha_1/\alpha_n)^k,\ldots,(\alpha_{n-1}/\alpha_n)^k\right)=
(x_1,\ldots,x_{n-1})$ is an arithmetic progression. Consequently, the set of
$k\in\Z$ such that no proper subsum of the left-hand side of (\ref{5.2})
vanishes, is the union of at most $A(n-1,1)$ arithmetic progressions.

Let $I$ be a proper, nonempty subset of $\{1,\ldots,n\}$ and consider those
solutions $k\in\Z$ of (\ref{5.2}) for which
\be
\label{5.5}
    \sum_{i\in I}a_i\,\alpha_i^k=0.
\ee
Each such $k$ also satisfies
\be
\label{5.6}
    \sum_{i\notin I}a_i\,\alpha_i^k=0.
\ee
Suppose $I$ has cardinality $\ell$. Since $a_i\neq 0$ $\,(i=1,\ldots,n)$, we
get $2\leq\ell\leq n-2$. By induction, the set of $k\in\Z$ with (\ref{5.5})
is the union of at most $W(\ell)$ arithmetic progressions. Also by induction,
the set of $k\in\Z$ with (\ref{5.6}) is the union of at most $W(n-\ell)$
arithmetic progressions. The intersection of two arithmetic progressions is
either empty, or again an arithmetic progression. In view of (\ref{5.6}) and
since $\ell\leq n-2$, $\,n-\ell\leq n-2$, the set of $k\in\Z$ with (\ref{5.2}),
(\ref{5.5}) is the union of at most
\begin{eqnarray*}
  W(\ell)W(n-\ell)&\leq&\exp\left((6\ell)^{3\ell}\right)\exp\left((6(n-\ell))^{3(n-\ell)}\right)\\
                  &\leq&\exp\left((6(n-1))^{3(n-1)}\right)\ \,=\ \,W(n-1)
\end{eqnarray*}
arithmetic progressions.

Taking into account all possible subsets $I$ of $\{1,\ldots,n\}$, we infer
that the set of solutions $k\in\Z$ of (\ref{5.2}) for which some subsum of
the left-hand side of (\ref{5.2}) vanishes is contained in the union of at
most $2^nW(n-1)$ arithmetic progressions. Recall from (\ref{5.4}) that the
set of $k\in\Z$ with (\ref{5.2}) for which no subsum of (\ref{5.2}) vanishes
is the union of at most $A(n-1,1)$ arithmetic progressions. So altogether, by
(\ref{5.1}) and  (\ref{5.4}), the set of solutions of (\ref{5.2}) is the union of
at most
$$
 \exp\left((6(n-1))^{3(n-1)}2\right)+2^n\exp\left((6(n-1))^{3(n-1)}\right)\leq\exp\left((6n)^{3n}\right)=W(n)
$$
arithmetic progressions.

Now suppose that $\{u_m\}$ is nondegenerate. Assume that ${\cal S}(u_m)$,
that is the set of solutions of (\ref{5.2}), contains an arithmetic
progression $\{a+vt\ |\ t\in\Z\}$ with $v\neq 0$. Then
\be
\label{5.7}
   a_1\,\alpha_1^a(\alpha_1^v)^t+\cdots+a_n\,\alpha_n^a(\alpha_n^v)^t=0
     \ \ \mbox{for every}\ \ t\in\Z.
\ee
Applying (\ref{5.7}) with $t=0,\ldots,n-1$, and observing that
$a_i\,\alpha_i^a\neq 0$ for $i=1,\ldots,n$, we infer that the Vandermonde
determinant $\,\mbox{det}\left(\alpha_i^{vt}\right)_{i=1,\ldots,n;\,t=0,\ldots,
n-1}$ is zero. This is possible only if there are $i\neq j$ with $\alpha_i^v=
\alpha_j^v$. But this contradicts the assumption that $\{u_m\}$ is
nondegenerate.

We conclude that ${\cal S}(u_m)$ does not contain an infinite arithmetic
progression. It follows that for nondegenerate $\{u_m\}$ the set ${\cal S}
(u_m)$ has cardinality $\leq W(n)$. This proves Theorem~1.2.\pagebreak

\section{A reduction}

We now turn to the proof of Theorem 2.1. Similarly as in the  argument used in
Section 3 in the deduction of Theorem 1.1, we claim that in order to prove 
Theorem 2.1 it
will suffice to show that any {\it{finite}} set $M$ of points $\by\in\overline
{\Q}^{\,n}$ satisfying (\ref{2.6}) and (\ref{2.10}) is contained in the union
of not more than
\be
\label{6.1}
    \exp\left((5n)^{3n}(r+1)\right)
\ee
proper linear subspaces of $\overline{\Q}^{\,n}$.

To verify this claim we prove:

\bl
Let $n\geq 2$ and $w\geq 1$ be integers{\rm .} Let $K$ be a field{\rm .} Let $N$ be a
subset of $K^n$ having the following property\/{\rm :}
\begin{center}
Any finite subset $M$ of $N$ is contained in the union\\ of not more than $w$
proper linear subspaces of $K^n${\rm .}\end{center}
Then $N$ itself is contained in the union of not more than $w$ proper linear
subspaces of $K^n${\rm .}
\el 

{\it Proof}.  By a subspace we shall mean a {\it{proper}} linear
subspace of $K^n$. 
Given a finite subset $M$ of $N$, we denote by $a(M)$ the minimum of the
quantities
$$
   \sum_{i=1}^w\dim T_i,
$$
where $\{T_1,\ldots,T_w\}$ runs through the collection of unordered $w$-tuples
of subspaces with
$$
   M\subset T_1\cup\ldots\cup T_w.
$$
Let ${\frak S}(M)$ be the collection of all $w$-tuples of subspaces
$\{T_1,\ldots,T_w\}$ with $M\subset T_1\cup\ldots\cup T_w$ and $\sum
\limits_{i=1}^w\dim T_i=a(M)$.

Now suppose $\{T_1,\ldots,T_w\}\in{\frak S}(M)$. Then $T_i$ $\,(i=1,\ldots,w)$
is generated by a subset of $M$. Otherwise we could replace $T_i$ by the
smaller subspace, generated by $T_i\cap M$, thus making $\sum\limits_{i=1}^w
\dim T_i$ smaller without affecting $M\subset T_1\cup\ldots\cup T_w$.

We may conclude that for each of the subspaces $T_i$ there are only finitely
many possibilities. Consequently, ${\frak S}(M)$ is finite. Denote its
cardinality by $b(M)$.

For any finite subset $M$ of $N$ we have $a(M)\leq(n-1)w$. Hence there is
such a subset $M$ for which $a(M)$ attains its maximum $a_0$, say. We now
choose among all finite subsets $M$ of $N$ having $a(M)=a_0$ a set $M_0$
such that
$$
    b(M_0)=\min_M\{b(M)\ |\ a(M)=a_0\}.
$$
If $M$ is any finite subset of $N$ with $M\supseteq M_0$ then ${\frak S}(M)=
{\frak S}(M_0)$. Indeed suppose $\{T_1,\ldots,T_w\}\in{\frak S}(M)$. So in
particular $M_0\subset T_1\cup\ldots\cup T_w$. On the other hand $\sum
\limits_{i=1}^w\dim T_i=a(M)\leq a_0$. The definition of $a_0$ implies that
$a(M)=a_0$ and therefore $\{T_1,\ldots,T_w\}\in{\frak S}(M_0)$; hence
${\frak S}(M)\subset{\frak S}(M_0)$. The inclusion ${\frak S}(M_0)\subset
{\frak S}(M)$ follows from the minimality of $b(M_0)$.

Pick $\{T_1,\ldots,T_w\}\in{\frak S}(M_0)$. We claim that
$$
    N\subset T_1\cup\ldots\cup T_w.
$$
Indeed let $\by\in N$ and consider the finite set $M=M_0\cup\{\by\}$. We have
shown that $\{T_1,\ldots,T_w\}\in{\frak S}(M)$. So in particular we have
$$
    \by\in T_1\cup\ldots\cup T_w.
$$
This proves our claim and the assertion of the lemma follows.

We now consider a finite set $M=\{\by_1,\ldots,\by_m\}$ of points $\by_i\in
\overline{\Q}^{\,n}$ satisfying (\ref{2.6}) and (\ref{2.10}). So we have
$$
  \by_i=\bx_i*\bz_i
$$
with
$$
  \bx_i\in\Gamma,\ \,\bz_i\in{(\,\overline{\Q}^{\,*})}^n,\ \,h(\bz_i)\leq
    n^{-1}\exp\left(-(4n)^{3n}\right)(1+h(\bx_i))\quad (i=1,\ldots,m).
$$
Let $F$ be a number field such that
$$
   \bx_i,\,\bz_i\in F^n\quad (i=1,\ldots,m).
$$
Write $\Gamma'$ for the subgroup of $\Gamma$ generated by $\bx_1,\ldots,\bx_m$.
Then $\Gamma'$ is a finitely generated subgroup of $(F^*)^n$ of rank $\leq r$.
Therefore, in order to prove (\ref{6.1}) for a finite set $M$ (and therefore
also Theorem 2.1) it will suffice to prove:

\proclaim{Proposition}
 Suppose $n\geq 2${\rm .} Let $F$ be a number field. Let $\Gamma$ be a finitely
generated subgroup of $(F^*)^n$ of rank $r${\rm .} Then the set of points $\by=
(y_1,\ldots,y_n)$ satisfying
\be
\label{6.2}
    y_1+\cdots+y_n=1,
\ee
\be
\label{6.3}
    \by=\bx*\bz\ \,\mbox{with}\ \,\bx\in\Gamma,\ \bz\in(F^*)^n,\ h(\bz)\leq
        n^{-1}\exp\left(-(4n)^{3n}\right)(1+h(\bx))
\ee
is contained in the union of not more than
\be
\label{6.4}
    \exp\left((5n)^{3n}(r+1)\right)
\ee
proper linear subspaces of $F^n${\rm .}
\endproclaim

The remainder of the paper deals with the proof of Proposition 6.2.
\pagebreak

\section{Heights in multiplicative groups}

For points $\ \bx=(x_1,\ldots,x_n)\in\overline{\Q}^{\,n}\setminus\{\b0\}\ $
and a number field $F$ such that $\ x_i\in F$ 
$(i=1,\ldots,n)$, we have defined in (\ref{2.3}) and (\ref{2.4}) respectively
the absolute multiplicative height
\be
\label{7.1}
    H(\bx)=\prod_{v\in M(F)}\max\{1,\|x_1\|_v,\ldots,\|x_n\|_v\}
\ee
as well as the absolute logarithmic height $h(\bx)=\log H(\bx)$. Thus, for
$\bx\in F^n$
\be
\label{7.2}
    h(\bx)=\sum_{v\in M(F)}\max\{0,\log\|x_1\|_v,\ldots,\log\|x_n\|_v\}.
\ee
Both, $H(\bx)$ as well as $h(\bx)$, do not depend upon the particular number
field $F$ such that $\bx\in F^n$. In the special case when $n=1$, (\ref{7.2})
yields for $x\in\overline{\Q}^{\,*}$ and a number field $F$ such that $x\in F$
\be
\label{7.3}
    h(x)=\sum_{v\in M(F)}\max\{0,\log\|x\|_v\}=\frac12\sum_{v\in M(F)}|\log
\|x\|_v|
\ee
(the last equation is a consequence of the product formula).

Then
$$
    h\left(\frac1x\right)=h(x),\ \ h(xy)\leq h(x)+h(y).
$$
For $\bx=(x_1,\ldots,x_n)\in{(\,\overline{\Q}^{\,*})}^n$ we define moreover
\be
\label{7.4}
    h_s(\bx)=\sum_{i=1}^nh(x_i).
\ee
Using (\ref{7.2})--(\ref{7.4}) we see that
\be
\label{7.5}
    h(\bx)\leq h_s(\bx)\leq n\,h(\bx).
\ee
Denoting as before by $*$ the product operation in ${(\overline{\Q}^{\,*})}^n$,
so that $(x_1,\ldots,x_n)*(y_1,\ldots,y_n)=(x_1y_1,\ldots,x_ny_n)$, we have
\be
\label{7.6}
   h(\bx*\by)\leq h(\bx)+h(\by)
\ee
and similarly for $h_s$. Further $h_s$ (but not $h$) is invariant under
replacement of  $\bx$ by its inverse $\bx^{-1}$ in ${(\,\overline{\Q}^{\,*})}^n$,
so that
\be
\label{7.7}
   h_s(\bx^{-1})=h_s(\bx).
\ee
From now on we fix the number field $F$. We let $\Gamma\subseteq(F^*)^n$ be a
finitely generated group of rank $r>0$. Let $\ba_1,\ldots,\ba_r$ be a set of
generators of $\Gamma$, so that the elements of $\Gamma$ are of the shape
\be
\label{7.8}
   \bx= \bxi*\ba_1^{u_1}*\ldots*\ba_r^{u_r}
\ee
where $(u_1,\ldots,u_r)$ runs through $\Z^r$, and $\bxi$ runs through the
torsion group $T(\Gamma)=\Gamma\cap U^n$ of $\Gamma$, where $U$ is the group
of roots of unity of $F$. For $\bu=(u_1\ldots,u_r)\in\Z^r$ set
\be
\label{7.9}
    \psi(\bu)=h_s(\ba_1^{u_1}*\ldots*\ba_r^{u_r}).
\ee
For $v\in M(F)$ put
$$
    \alpha_{ijv}=\log\|a_{ij}\|_v\quad (1\leq i\leq r,\ 1\leq j\leq n)
$$
where $\ \ba_i=(a_{i1},\ldots,a_{in})$. Then by the product formula
$\ \sum\limits_{v\in M(F)}\alpha_{ijv}=0$ $\ (1\leq i\leq r$, 
$1\leq j\leq n)$.
Let $S$ be the subset of $M(F)$ consisting of those $v$'s such that
$\alpha_{ijv}\neq 0$ for some pair $i,\,j$ $\,(1\leq i\leq r,\ 1\leq j\leq n)$.
Then also $\sum\limits_{v\in S}\alpha_{ijv}=0$ $\,(1\leq i\leq r,\ 1\leq j
\leq n)$. For $\bxi\in\R^r$ put
\be
\label{7.10}
    g_{jv}(\bxi)=\sum_{i=1}^r\alpha_{ijv}\,\xi_i\quad (1\leq j\leq n,\ v\in
        M(F));
\ee
then again
\be
\label{7.11}
   \sum_{v\in S}g_{jv}(\bxi)=0\ \ (1\leq j\leq n)\ \ \mbox{and}\ \ g_{jv}
   (\bxi)=0\ \ \mbox{for}\ \ v\notin S\ \ (1\leq j\leq n).
\ee
Since by (\ref{7.8})
\be
\label{7.12}
    \log\|x_j\|_v=\log\|a_{1j}^{u_1}\ldots a_{rj}^{u_r}\|_v=\sum_{i=1}^r
      \alpha_{ijv}\,u_i=g_{jv}(\bu), 
\ee   
we have from (\ref{7.3}), (\ref{7.9}),
\begin{eqnarray*}
 \psi(\bu)\ \,=\ \,\sum^n_{j=1}h(a_{1j}^{u_1}\ldots a_{rj}^{u_r})&=&\frac12
  \,\sum^n_{j=1}\ \sum_{v\in M(F)}\,|g_{jv}(\bu)|\\
             &=&\frac12\,\sum^n_{j=1}\ \sum_{v\in S}\,|g_{jv}(\bu)|.
\end{eqnarray*}
More generally, for $\bxi\in\R^r$ set
\be
\label{7.13}
  \psi(\bxi)\,=\,\frac12\,\sum_{v\in M(F)}\ \sum^n_{j=1}\,|g_{jv}(\bxi)|.
\ee
Then
\begin{itemize}
 \item[(a)] $\ \psi(\bxi)\geq 0\ \ \mbox{for}\ \ \bxi\in\R^r$,
 \item[(b)] $\ \psi(\alpha\,\bxi)=|\alpha|\,\psi(\bxi)\ \ \mbox{for}\ \ \bxi
                 \in\R^r,\ \alpha\in\R$,
 \item[(c)] $\ \psi(\bxi+\bet)\leq\psi(\bxi)+\psi(\bet)\ \ \mbox{for}\ \ \bxi,
                 \bet\in\R^r$.
\end{itemize}
Since $\ba_1,\ldots,\ba_r$ are multiplicatively independent, the components
of $\ba_1^{u_1}*\ldots*\ba_r^{u_r}$ will all be roots of unity only if $\bu
=\b0$. Therefore according to Dobrowolski \cite{6}, for $\bu\in\Z^r\setminus
\{\b0\}$  we have
$$
  \psi(\bu)>\frac{c_1}d\,(\log\log 3d\,/\log 3d)^3
$$
where $d=[F:\Q\,]$ and where $c_1>0$ is an absolute constant. In particular,
there is a constant $c>0$ such that
\begin{itemize}
 \item[(d)] $\ \psi(\bu)\geq c>0\ \ \mbox{for}\ \ \bu\in\Z^r\setminus\{\b0\}$.
\end{itemize}
In \cite{25}, Lemma 3 it is shown that since the function $\psi$ satisfies
(a)--(d), the set $\Psi\subset\R^r$ given by
\be
\label{7.14}
    \Psi=\{\bxi\in\R^r\,|\ \psi(\bxi)\leq 1\}
\ee
is a {\it{symmetric{\rm ,} convex body}}.

\section{Special points}

Let $F$, $\,\Gamma$, $\,\ba_1,\ldots,\ba_r\,$ be as in Section 7. When $\bx
\in\Gamma$, set
\be
\label{8.1}
    h=h(\bx),\ \ H=H(\bx)=e^h,\ \ h_s=h_s(\bx).
\ee
Express $\bx$ as in (\ref{7.8}). So if $\bx\in\Gamma$ and $\bu\in\Z^r$ are
related by (\ref{7.8}), we have (\ref{7.9}), i.e.,
\be
\label{8.2}
    h_s=h_s(\bx)=\psi(\bu).
\ee
Let $\Psi=\{\bxi\in\R^r\,|\ \psi(\bxi)\leq 1\}$ be the set (\ref{7.14}). Put
\be
\label{8.3}
    q=4n.
\ee
{\it{Given $\br\in\R^r${\rm ,} an element $\bx\in\Gamma$ will be called $\br$-special
if}} $h>0$ ({\it{in}} (\ref{8.1})) {\it{and if}}
\be
\label{8.4}
    \bu\in\frac hq\,\Psi+h\,\br.
\ee
The right-hand side of (\ref{8.4}) signifies $\,\displaystyle{\frac hq}\,\Psi$
translated by $h\,\br$.

We quote Lemma 8.1 of \cite{23}. 

\bl
Let $\Phi$ be a symmetric convex body in $\R^r$. Suppose $\lambda>0${\rm .} Then
$\lambda\,\Phi$ can be covered by not more than
\be
\label{8.5}
    (2\lambda+4)^r
\ee
translates of $\Phi${\rm .}
\el

We apply Lemma 8.1 with $\Phi$ replaced by $\,\displaystyle{\frac1q}\,\Psi$
and with $\lambda
\,\Phi$ replaced by $n\,\Psi$. We may conclude that $n\,\Psi$ may be covered by not more than
$$
   (2qn+4)^r=Z
$$
translates of $\,\displaystyle{\frac1q}\,\Psi$, say by $\,\displaystyle
{\frac1q}\,\Psi+\br_i$ $\,(i=1,\ldots,Z)$.

Now when $\bx$ satisfies (\ref{8.1}), then by (\ref{7.5}), (\ref{7.9}),
(\ref{8.2}) the point $\bu\in\Z^r$ related to $\bx$ via (\ref{7.8}) lies in
$$
    h_s\,\Psi\subset h\,n\,\Psi.
$$
Thus $\bx$ is special for at least one of $\br_1,\ldots,\br_Z$. We have shown:

\proclaim{{C}orollary}
 There exist elements $\br_1,\ldots,\br_Z\in \R^r$ with 
\be
\label{8.6}
    Z=(2qn+4)^r
\ee
 such that any $\bx\in\Gamma$ is special for at least one of $\br_1,\ldots,
\br_Z${\rm .}
\endproclaim

We remark moreover that our construction implies that we may take
$\br_1,\ldots,\br_Z$ with
\be
\label{8.7}
  \br_i\in\left(n+\frac1q\right)\!\Psi\quad (i=1,\ldots,Z).
\ee
In the sequel, we will apply the material developed so far to the solutions
$\by=\bx*\bz$ of (\ref{6.2}), (\ref{6.3}).

\section{Properties of large special solutions}
 
We now study solutions $\by$ of (\ref{6.2}), (\ref{6.3}).
\vglue4pt
{\it{A solution $\by$ will be called large if it has a representation $\by=
\bx*\bz$ as in}} (\ref{6.3}) {\it{such that}}
\be
\label{9.1}
    h(\bx)>4n\log n.
\ee
{\it{Solutions $\by$ of}} (\ref{6.2}), (\ref{6.3}) {\it{that are not large
will be called small.}}
\vglue4pt
If the group $\Gamma$ has rank $0$, then all elements $\bx\in\Gamma$ have
$h(\bx)=0$. So, large solutions only exist when $\mbox{rank}\,\Gamma>0$. 
\vglue4pt
{\it{A solution $\by$ of}} (\ref{6.2}), (\ref{6.3}) {\it{is called
$\br$-special if{\rm ,} with $\bx,\bz$ as in}} (\ref{6.3}), {\it{the point $\bx$ is
$\br$-special.}}
\vglue4pt

In this section we derive properties of large $\br$-special solutions $\by$.
This will allow us in Section 10 to deduce an upper bound for the number of
subspaces needed to cover the set of large solutions of (\ref{6.2}),
(\ref{6.3}).

Suppose that $\br\in\left(n+\displaystyle{\frac1q}\right)\!\Psi$ is fixed. Set
\be
\label{9.2}
   m_{jv}=\left\{
       \begin{array}{ll}
                  g_{jv}(\br) & (v\in M(F),\ 1\leq j\leq n)\\
                  0           &  (v\in M(F),\ j = 0).
       \end{array}
   \right.
\ee
In view of (\ref{7.11}) we have
\be
\label{9.3}
   \sum_{v\in S}m_{jv}=0\quad(j=0,\ldots,n),\quad m_{jv}=0\ \,\mbox{for}
     \ \,v\notin S,\ j=0,\ldots,n. \hskip.25in
\ee
Further, since $\br\in\left(n+\displaystyle{\frac1q}\right)\!\Psi$, by
(\ref{9.2}),
(\ref{7.13}) and the
definition of $\Psi$ in (\ref{7.14}),
\be
\label{9.4}
    \sum_{v\in M(F)}\ \sum^n_{j=0}\,|m_{jv}|=2\psi(\br)\leq 2
     \left(n+\frac1q\right).
\ee
Now let $\bx\in\Gamma$ be $\br$-special, so that with $\bu$ as in (\ref{7.8})
we have (\ref{8.4}) with $h=h(\bx)$. Then for any $v\in M(F)$ and for $j=1,
\ldots,n$
\be
\label{9.5}
   g_{jv}(\bu)=h(g_{jv}(\br)+q^{-1}\,g_{jv}(\bxi))=h\,m_{jv}+\frac hq
      \,g_{jv}(\bxi)
\ee
with some $\bxi\in\Psi$. Writing $g_{0v} 
(\bxi)=0$ for $v\in M(F)$ and for $\bxi\in\R^r$, (\ref{9.5}) will be true for
$j=0$ as well.
 
It follows from (\ref{9.2}), (\ref{9.5}) that
\be
\label{9.6}
   \sum_{v\in M(F)}\ \sum^n_{j=0}\ |g_{jv}(\bu)-h\,m_{jv}|=\frac hq
      \,\sum_{v\in M(F)}\ \sum^n_{j=0}\ |g_{jv}(\bxi)|\leq\frac hq\,.
\ee
For $v\in M(F)$ let $L_0^{(v)},\ldots,L_n^{(v)}$ be the linear forms in $\bY=
(Y_1,\ldots,Y_n)$ given by
\begin{eqnarray}
\label{9.7}
     L^{(v)}_0(\bY) & = & Y_1 + \ldots + Y_n,\\
     L^{(v)}_1(\bY) & = & Y_1,\nonumber \\
     \vdots          & & \nonumber\\
     L^{(v)}_n(\bY) & = & Y_n\,. \nonumber
\end{eqnarray}

\bl
Let $\br$ be as above{\rm .} There are $n$\/{\rm -}\/element subsets $\Is(v)$ of
$\{0,\ldots,n\}$ defined for $v\in M(F)$ and there are numbers $\ell_{jv}$
$\,(v\in M(F),\ j\in\Is(v))$ with the following properties{\rm .}
\begin{eqnarray}
\label{9.8}
   \Is(v)=\{1,\ldots,n\}\phantom{,} & & \mbox{for}\quad v\notin S,\\
\label{9.9}
   \ell_{jv}=0\phantom{,} & & \mbox{for}\quad v\notin S,\ \ j\in\Is(v),\\
\label{9.10}
   \sum_{v\in M(F)}\ \sum_{j\in\Is(v)}\,\ell_{jv}=0, &  &
      \sum_{v\in M(F)}\ \sum_{j\in\Is(v)}\,|\ell_{jv}|\leq 1. 
\end{eqnarray}
Moreover{\rm ,} any large $\br$\/{\rm -}\/special solution $\by$ of {\rm (\ref{6.2}),
(\ref{6.3})} satisfies the inequality
\be
\label{9.11}
   \prod_{v\in M(F)}\ \max_{j\in\Is(v)}\,\left\{\frac{\|L_j^{(v)}(\by)\|_v}
      {Q^{\ell_{jv}}}\right\}\,\leq\,Q^{-\frac1{2n(4n+1)}},
\ee
where $Q=H(\bx)^{4n+1}${\rm .} Here $\bx\in\Gamma$ is a point in the representation
$\by=\bx*\bz$ according to {\rm (\ref{6.3}).}
\el

{\it Proof}.  For $v\in S$ let $j(v)\in\{0,\ldots,n\}$ be a subscript with
\be
\label{9.12}
    m_{j(v),v}=\max\,\{m_{0v},\ldots,m_{nv}\}.
\ee
We define $\Is(v)=\{0,\ldots,n\}\setminus\{j(v)\}$ $\,(v\in S)$. For
$v\notin S$, $\,\Is(v)\,$ is already defined in (\ref{9.8}). By our definition
of $S$ in Section 7, any solution $\by=\bx*\bz$ as in (\ref{6.3}) has by
(\ref{7.11}), (\ref{7.12})
\begin{eqnarray}
\label{9.13}
  h\ \ =\ \ h(\bx)&=&\sum_{v\in S}\ \max\,\{0,\,\log\|x_1\|_v,\ldots,
   \log\|x_n\|_v\}\\ \nonumber
  &=&\sum_{v\in S}\ \max\,\{g_{0v}(\bu),\,g_{1v}(\bu),\ldots,g_{nv}(\bu)\}.
\end{eqnarray}
Given $\bx$, pick for each $v\in S$ an element $i(v)\in\{0,\ldots,n\}$ with
$g_{i(v),v}(\bu)=\max\,\{g_{0v}(\bu),\ldots,g_{nv}(\bu)\}$. Then by
(\ref{9.13})
$$
   \sum_{v\in S}g_{i(v),v}(\bu)=h.
$$ 
Thus in view of (\ref{9.6}), we may infer that
$$
   h\,\sum_{v\in S}m_{i(v),v}\,\geq\,\sum_{v\in S}g_{i(v),v}(\bu)-\frac hq=
     h\left(1-\frac1q\right).
$$
In particular, by (\ref{9.12}) we obtain
\be
\label{9.14}
   \sum_{v\in S}m_{j(v),v}\,\geq\, 1-\frac1q\,.
\ee
Let $s$ be the cardinality of $S$ and write
\be
\label{9.15}
   \gamma=\frac1{ns}\,\sum_{v\in S}m_{j(v),v}.
\ee
We now define numbers $c_{jv}$ $\,(v\in M(F),\ j\in\Is(v))$ by
\be
\label{9.16}
   c_{jv}=\left\{
   \begin{array}{ll}
          m_{jv}+\gamma & \mbox{for}\quad v\in S,\ \, j\in\Is(v)\\
          0             & \mbox{for}\quad v\notin S,\ \, j\in\Is(v).
   \end{array}\right.
\ee
We infer from (\ref{9.3}), (\ref{9.4}), (\ref{9.15}) that
\be
\label{9.17}
   \sum_{v\in M(F)}\ \sum_{j\in\Is(v)}c_{jv} = 0,\qquad
   \sum_{v\in M(F)}\ \sum_{j\in\Is(v)}|c_{jv}|\leq 4\left(n+\frac1q\right).\hskip.25in
\ee
So far we have only used the fact that our solution $\by=\bx*\bz$ of
(\ref{6.2}), (\ref{6.3}) is $\br$-special.

However, $\by$ is also supposed to be large. Under this additional hypothesis
we now derive an upper bound for the quantity
$$
   A=\prod_{v\in M(F)}\max_{j\in\Is(v)}\left\{\frac{\|L_j^{(v)}(\by)\|_v}
       {H^{c_{jv}}}\right\},
$$
where $H=H(\bx)$.

Write $\by=(y_1,\ldots,y_n)$, $\,\bx=(x_1,\ldots,x_n)$, $\,\bz=(z_1,\ldots,
z_n)$. Put $y_0=x_0=z_0=1$. Notice that by (\ref{6.2}) and (\ref{9.7}) we
then have for each $v\in M(F)$
$$
   L_j^{(v)}(\by)=y_j\ \ \mbox{for}\ \ j=0,\ldots,n.
$$
Hence by (\ref{9.14}), (\ref{9.15}), (\ref{9.16}),
\begin{eqnarray}
\label{9.18}\qquad
  A&=&\prod_{v\in M(F)}\max_{j\in\Is(v)}\left\{\frac{\|y_j\|_v}{H^{c_{jv}}}
                                                   \right\} \\
   &=&H^{-\frac1n\sum_{v\in S}m_{j(v),v}}\prod_{v\in M(F)}\max_{j\in\Is(v)}
         \left\{\frac{\|y_j\|_v}{H^{m_{jv}}}\right\}\nonumber\\
   &\leq&H^{-\frac1n+\frac1{nq}}\prod_{v\in M(F)}\max_{j\in\Is(v)}
         \left\{\frac{\|y_j\|_v}{H^{m_{jv}}}\right\}\nonumber\\
   &\leq&H^{-\frac1n+\frac1{nq}}\left(\prod_{v\in M(F)}\max_{0\leq j\leq n}
         \|z_j\|_v\right)\left(\prod_{v\in M(F)}\max_{0\leq j\leq n}
         \left\{\frac{\|x_j\|_v}{H^{m_{jv}}}\right\}\right)\nonumber\\
   &=&H^{-\frac1n+\frac1{nq}}H(\bz)\prod_{v\in M(F)}\max_{0\leq j\leq n}
         \left\{\frac{\|x_j\|_v}{H^{m_{jv}}}\right\}.\nonumber
\end{eqnarray}
 (\ref{6.3}) and (\ref{9.1}) entail
\begin{eqnarray}
\label{9.19}
&&\\
  H(\bz)&\leq&\exp\left(n^{-1}\exp\left(-(4n)^{3n}\right)(1+h(\bx))\right) \nonumber\\ \nonumber
         &\leq&\exp\left(n^{-1}\exp\left(-(4n)^{3n}\right)\left((4n\log n)^{-1}+1\right)h(\bx)\right)
               \ \,\leq\ \,H^{1/(8n)}.
\end{eqnarray}
On the other hand by (\ref{7.12}) and (\ref{9.2}), (\ref{9.5}), (\ref{9.6}),
\be
\label{9.20}
   \prod_{v\in M(F)}\max_{0\leq j\leq n}\left\{\frac{\|x_j\|_v}{H^{m_{jv}}}
      \right\}\,\leq\,H^{1/q}.
\ee
Combination of (\ref{9.18})--(\ref{9.20}) yields, with our value $q$ from
(\ref{8.3}),
\be
\label{9.21}
   \prod_{v\in M(F)}\max_{j\in\Is(v)}\left\{\frac{\|L_j^{(v)}(\by)\|_v}
     {H^{c_{jv}}}\right\}\,\leq\,H^{-\frac1n+\frac1{nq}+\frac1{8n}+\frac1q}
     \,\leq\,H^{-\frac1{2n}}. \hskip.5in
\ee
We now renormalize with our parameter $Q=H(\bx)^{4n+1}=H^{4n+1}$. Writing
$$
   \ell_{jv}\,=\,\frac{c_{jv}}{4n+1}\qquad(v\in M(F),\ \,j\in\Is(v))
$$
we obtain, with (\ref{9.16}), (\ref{9.17}) and with $q$ as in (\ref{8.3}),
assertions (\ref{9.9}) and (\ref{9.10}). Moreover, (\ref{9.21}) gives
(\ref{9.11}).

\section{Large solutions}

To deal with the large solutions, we use the absolute version of the
Subspace Theorem, due to Evertse and Schlickewei \cite{12}. The following
Proposition~10.1 is a very special case of Theorem 2.1 of \cite{12}.

For $v\in M(F)$ let the linear forms $L_0^{(v)}(\bY),\ldots,L_n^{(v)}(\bY)$
be as in (\ref{9.7}). Moreover, let $\Is(v)$ and the tuple $(\ell_{jv})$
$\,(v\in M(F),\,j\in\Is(v))$ be as in Lemma~9.1.

\bp
Suppose $0<\delta<1${\rm .} There are proper linear subspaces $T_1,\ldots,T_t$
of $F^n$ with
\be
\label{10.1}
   t\leq 2^{2(n+9)^2}\delta^{-n-4}
\ee
with the following property\/{\rm :} 
As $Q$ runs through the values satisfying
\be
\label{10.2}
   Q>n^{2/\delta},
\ee
the set of solutions $\by\in F^n$ of the inequalities
\be
\label{10.3}
  \prod_{v\in M(F)}\max_{j\in\Is(v)}\left\{\frac{\|L^{(v)}_j(\by)\|_v}
     {Q^{\ell_{jv}}}\right\}\leq\,Q^{-\delta}
\ee
is contained in the union
$$
   T_1\cup\ldots\cup T_t.
$$
\ep

We apply Proposition 10.1 with $Q=H(\bx)^{4n+1}$ (where $\by=\bx*\bz$ with
$\bx\in\Gamma$ according to (\ref{6.3})), and with $\delta=
\displaystyle{\frac1{2n(4n+1)}}\,$.
By Lemma 9.1, given $\br$, any large $\br$-special solution $\by$ of
(\ref{6.2}), (\ref{6.3}) satisfies (\ref{10.3}) with sets $\Is(v)$ and a tuple
$(\ell_{jv})$ $\,(v\in M(F),\,j\in\Is(v))$ which depend only on $\br$.

With our values of $Q$ and $\delta$, (\ref{10.2}) becomes $H(\bx)^{4n+1}>
n^{4n(4n+1)}$, or equivalently
$$
   h(\bx)>4n\log n.
$$
In view of (\ref{9.1}) this means that Proposition 10.1 is adequate to deal
with the large $\br$-special solutions $\by$ of (\ref{6.2}), (\ref{6.3}).
\pagebreak

By (\ref{10.1}), a single $\br$ gives rise to not more than
$$
 2^{2(n+9)^2}(8n^2+2n)^{n+4}
$$
subspaces. Using Corollary 8.2 and the definition of $q$ in (\ref{8.3}) we
obtain

\proclaim{{C}orollary}  The set of large solutions $\by$ of  {\rm (\ref{6.2}), (\ref{6.3})}  is
contained in the union of not more than 
$$
   2^{2(n+9)^2}(8n^2+2n)^{n+4+r}
$$
 proper linear subspaces of $F^n${\rm .}
\endproclaim

\section{Small solutions}

We still have to deal with the small solutions $\by$ of (\ref{6.2}),
(\ref{6.3}).
For this purpose we use results on the number of points on varieties which
have small height. The first {\it{explicit}} estimate in that context is due to
W. Schmidt \cite{26}. We quote here a special case of Theorem 4 of \cite{26}.

\bp
Let $\bb=(b_1,\ldots,b_n)\in{(\,\overline{\Q}^{\,*})}^n${\rm .} Put
\be
\label{11.1}
    q_0(n)=\exp\left((4n)^{3n}\right).
\ee
Then the equation
\be
\label{11.2}
    b_1\,w_1+\cdots+b_n\,w_n=1
\ee
has at most $q_0(n)$ nondegenerate solutions $\bw=(w_1,\ldots,w_n)\in
{(\,\overline{\Q}^{\,*})}^n$ with
\be
\label{11.3}
    h_s(\bw)<q_0(n)^{-1}.
\ee
\ep

We remark that S. David and P. Philippon \cite{4}, \cite{5} recently have
proved a sharpening of Proposition 11.1. They have shown that with
$$
   q_1(n)=2^{2^{c_1n}},
$$
where $c_1$ is an explicit absolute constant, equation (\ref{11.2}) has at
most $q_1(n)$ nondegenerate solutions $\bw$ with
$$
   h_s(\bw)<q_1(n)^{-3/4}.
$$
Here, we will give details on the basis of Proposition 11.1.

W. Schmidt, in Theorem 5 of \cite{26}, also has derived an upper bound for the
number of nondegenerate solutions $\bw$ of (\ref{11.2}) when $\bw$ lies in a
group $\Gamma\subset{(\,\overline{\Q}^{\,*})}^n$ of rank $r$ and has
\be
\label{11.4}
    h_s(\bw)\leq C.
\ee
In our context we ask for the number of points $\by=\bx*\bz$ satisfying
(\ref{6.2}), (\ref{6.3}). But our $\by$ only ``essentially'' belongs to
$\Gamma$ (in the sense defined by (\ref{6.3})). Moreover, instead of
(\ref{11.4}), which in our context would be $h_s(\by)\leq C$, we only have a
weaker hypothesis of type
\be
\label{11.5}
    h_s(\bx)\leq C.
\ee
To derive a bound in this more general setting, we follow the argument given
in \cite{26}.

By (\ref{6.3}), $\by=\bx*\bz$ with
\be
\label{11.6}
   \bx\in\Gamma\ \ \mbox{and}\ \ h(\bz)\leq n^{-1}\exp\left(-(4n)^{3n}\right)(1+h(\bx)).
\ee
 Suppose first that $\mbox{rank}\,\Gamma=0$. 
Then $h(\bx)=0$. Therefore, by (\ref{7.5}), (\ref{7.6}),
\begin{eqnarray}
\label{11.7}
    h_s(\by)&\leq&h_s(\bx)+h_s(\bz)\ \,=\ \,h_s(\bz)\\
            &\leq&n\,n^{-1}\exp\left(-(4n)^{3n}\right)(1+h(\bx)) \nonumber\\
&=&
                  \exp\left(-(4n)^{3n}\right)\ \,=\ \,q_0(n)^{-1}.\nonumber
\end{eqnarray}
We apply Proposition 11.1 with $\bb=(1,\ldots,1)$ and conclude that (\ref{6.2})
does not have more than
\be
\label{11.8}
    q_0(n)
\ee
nondegenerate solutions $\by$ satisfying (\ref{11.7}). We point out that our
choice of the function $n^{-1}\exp\left(-(4n)^{3n}\right)$ in (\ref{6.3}) is motivated
uniquely to guarantee (\ref{11.7}).

 We now treat the case when $r=\mbox{rank}\,\Gamma>0$. 
Hypothesis (\ref{11.5}), in view of (\ref{7.5}) and (\ref{9.1}), now reads as
\be
\label{11.9}
    h_s(\bx)\leq 4n^2\log n.
\ee
Let $\bu\in\Z^r$ be the point related to $\bx\in\Gamma$ by (\ref{7.8}).
Combination of (\ref{7.9}) and (\ref{11.9}) gives
\be
\label{11.10}
    \psi(\bu)\leq 4n^2\log n.
\ee
We quote Lemma 4 of \cite{25}.

\proclaim{Lemma}
 Let $\psi:\R^r\rightarrow\R$ be a function satisfying  {\rm (a)--(d)}
 in Section~{\rm 7.}   Let $U$ be a set of points in $\R^r$ such that 
\be
\label{11.11}
    \psi(\bu-\bv)\geq\delta_0>0
\ee
 for $\bu\neq\bv$ in $U${\rm .} Then the number of $\bu\in U$ with 
\be
\label{11.12}
    \psi(\bu)\leq C
\ee
 is 
\be
\label{11.13}
    \leq ((2C/\delta_0)+1)^r.
\ee
\endproclaim

Let $V$ be the subset of points $\bu\in\Z^r$ satisfying (\ref{11.10}). We
apply Lemma~11.2 with $U$ being a maximal subset of $V$ such that
\be
\label{11.14}
    \psi(\bu-\bv)\geq\frac12\,q_0(n)^{-1}\ \ \mbox{for}\ \ \bu\neq\bv
       \ \ \mbox{in}\ \ U.
\ee
Here $q_0(n)$ is as in Proposition 11.1. So we take $C=4n^2\log n$ and
$\delta_0=\displaystyle{\frac12}\,q_0(n)^{-1}$. By (\ref{11.13}) and
(\ref{11.1}) we may
infer that $U$ has cardinality
\be
\label{11.15}
    |U|\leq\left(16n^2(\log n)q_0(n)+1\right)^r\leq\left(16n^3\,q_0(n)\right)^r.
\ee
Moreover by the definition of $U$, for any $\bu\in\Z^r$ satisfying
(\ref{11.10}), there exists $\bu_0\in U$ such that
\be
\label{11.16}
    \psi(\bu-\bu_0)<\frac12\,q_0(n)^{-1}.
\ee
Again using (\ref{7.8}), (\ref{7.9}) we may infer that there is a subset
$\Delta$ of $\Gamma$ with cardinality
\be
\label{11.17}
    |\Delta|\leq\left(16n^3\,q_0(n)\right)^r
\ee
such that for any $\bx\in\Gamma$ with (\ref{11.9}) there is an element
$\bb\in\Delta$ having
\be
\label{11.18}
    h_s\left(\bx*\bb^{-1}\right)<\frac12\,q_0(n)^{-1}.
\ee
Now let $\by$ be a small solution of (\ref{6.2}), (\ref{6.3}), i.e., a
solution with $h(\bx)\leq 4n\log n$.

We choose $\bb\in\Delta$ satisfying (\ref{11.18}). Combination of (\ref{6.3}),
(\ref{7.5}), (\ref{7.6}), (\ref{11.1}), (\ref{11.18}) yields
\begin{eqnarray}
\label{11.19}
&&\\
    h_s(\by*\bb^{-1})&\leq&h_s(\bz)+h_s(\bx*\bb^{-1})\ \,\leq\ \,n\,h(\bz)+
                                             h_s(\bx*\bb^{-1}) \nonumber\\ \nonumber
                     &\leq&n\exp\left(-(5n)^{3n}\right)(1+4n\log n)+\frac12\,q_0(n)^{-1}
                                             \ \,<\ \,q_0(n)^{-1}.
\end{eqnarray}
We conclude that for any small nondegenerate solution $\by$ of (\ref{6.2}),
(\ref{6.3}) there exists $\bb\in\Delta$ with (\ref{11.19}).

Write $\bw=\by*\bb^{-1}$. Then $\bw$ is a solution of (\ref{11.2}),
(\ref{11.3}). By Proposition~11.1, given $\bb$, there are at most $q_0(n)$
points $\bw$ with (\ref{11.2}), (\ref{11.3}). We may conclude that each
$\bb\in\Delta$ gives rise to at most $q_0(n)$ nondegenerate small solutions
$\by$ of (\ref{6.2}).

Introducing the factor $(16n^3\,q_0(n))^r$ from (\ref{11.17}) for the number
of possible choices of $\bb$, we see that altogether we cannot have more than
\be
\label{11.20}
    q_0(n)\left(16n^3\,q_0(n)\right)^r
\ee
nondegenerate small solutions.

Comparing (\ref{11.20}) with (\ref{11.8}) we observe that indeed the bound
(\ref{11.20}) is true for any value of $r=\mbox{rank}\,\Gamma$.

All other {\it{small}} solutions are degenerate, i.e., some subsum on the
left-hand side of (\ref{6.2}) vanishes. The number of subsums is $\leq 2^n$.
Hence the degenerate solutions may be covered by the union of $\leq 2^n$
proper linear subspaces.

To summarize, we have proved:

\proclaim{{C}orollary}  The set of small solutions of  {\rm (\ref{6.2}), (\ref{6.3})}   is
contained in the union of not more than 
$$
   2^n+q_0(n)\left(16n^3q_0(n)\right)^r
$$
 proper linear subspaces of $F^n$. Here $q_0(n)$ is given by  {\rm (\ref{11.1}),}
  i.e.{\rm ,} $q_0(n)=\exp\left((4n)^{3n}\right)${\rm .}
\endproclaim

\section{Proof of Proposition 6.2}

We collect the results of Sections 10 and 11. From Corollary 10.2 we get not more than
$$
     2^{2(n+9)^2}\left(8n^2+2n\right)^{n+4+r}
$$
subspaces for the large solutions. From Corollary 11.2 we obtain not more than
$$
  2^n+\left(16n^3\right)^r\left(\exp((4n)^{3n})\right)^{r+1}
$$
subspaces for the small solutions. 
Therefore, to cover the set of all solutions $\by$ of (\ref{6.2}), (\ref{6.3})
$$
 2^{2(n+9)^2}\left(8n^2+2n\right)^{n+4+r}+2^n+\left(16n^3\right)^r\left(\exp((4n)^{3n})\right)^{r+1}
    \,<\,\exp\left((5n)^{3n}(r+1)\right)
$$
subspaces will suffice. This completes the proof of Proposition 6.2.

\bye
  
%  \vspace{2em}
%\pagebreak    
    \noindent 
%     J.-H. Evertse\\
%     Dept. of Mathematics and Computer Science\\
%     Rijksuniversiteit te Leiden\\
%     Niels Bohrweg 1\\
%     P.O. Box 9512\\
%     2300 RA Leiden\\
%     The Netherlands\\
%     e-mail: evertse@math.leidenuniv.nl
%     
%       \vspace{1em}
%
%       \noindent 
%     H.P. Schlickewei\\
%     Fachbereich Mathematik\\
%     Universit\"at Marburg\\
%     Lahnberge\\
%     35032 Marburg \\
%     Germany\\
%     e-mail: hps@mathematik.uni-marburg.de\
%     
%       \vspace{1em}
%
%  \noindent      
%     W. M. Schmidt\\
%     Dept. of Mathematics\\
%     University of Colorado\\
%     Boulder, CO 80309-0395\\
%     U.S.A.\\
%     e-mail: schmidt@euclid.colorado.edu
\bigskip
\centerline{(Received August 8, 2000)}
\end{document}